\documentclass[12pt]{aptpub}
\usepackage{color}
 \usepackage{amsmath}
 \usepackage{amsfonts}
\authornames{O. Zhao, M. Woodroofe and D. Voln\' y}
\shorttitle{CLT for Reversible Processes}

\def\be{\begin{equation}}
\def\ee{\end{equation}}

\def\R{\mathbb R}
\def\Z{\mathbb Z}

\def\ed{\end{document}}

\def\<{\langle}
\def\>{\rangle}
\def\eqalign#1{\null\,\vcenter{\openup\jot\ialign
              {\strut\hfil$\displaystyle{##}$&$\displaystyle{{}##}$
               \hfil\crcr#1\crcr}}\,}

\begin{document}
\title{A Central Limit Theorem  For Reversible Processes  With Non-linear Growth of Variance\footnote{Corresponding author: Ou Zhao (ouzhao@stat.sc.edu)} }

\authorone[USC-Columbia]{Ou Zhao}
\addressone{Department of Statistics, University of South Carolina, 1523 Greene Street, Columbia, SC 29208, USA.} 

\authortwo[University of Michigan]{Michael Woodroofe}

\addresstwo{Department of Statistics and Mathematics, University of
 Michigan, 275 West Hall, 1085 South University, Ann Arbor, MI 48109, USA.  Email address: michaelw@umich.edu}

\authorthree[Universit\'e de Rouen]{Dalibor Voln\' y}

\addressthree{Laboratoire de math\' ematiques Rapha\" el Salem, UMR CNRS 6085, 
 Universit\' e de Rouen, France.  Email address: Dalibor.Volny@univ-rouen.fr}

\begin{abstract}

	Kipnis and Varadhan showed that for an additive functional, $S_n$ say, of a reversible Markov chain the condition $E(S_n^{2})/n \to \kappa \in (0,\infty)$ implies the convergence of the conditional distribution of $S_n/\sqrt{E(S_n^{2}})$, given the starting point, to the standard normal distribution.  We revisit this question under the weaker condition, $E(S_n^{2}) = n\ell(n)$, where $\ell$ is a slowly varying function.   It is shown by example that the conditional distribution of $S_n/\sqrt{E(S_n^{2}})$ need not converge to the standard normal distribution in this case; and sufficient conditions for convergence to a (possibly non-standard) normal distribution are developed.

\end{abstract}

	
% {\it Key words and phrases}: conditional distributions; Markov chains; Metropolis Hastings; self-adjoint operators; slowly varying functions.

% \medskip\noindent
% {\it 2000 AMS Classification}: Primary 60F05; Secondary 60J05.

\bigskip
\section{Introduction}\label{sect:intro}

	Consider a reversible Markov chain $\ldots W_{-1},W_0,W_1,\ldots$, defined on a probability space $(\Omega,{\cal A},P)$, with a Polish state space ${\cal W}$, transition function $Q$, and marginal distribution $\pi$.  Thus, $\pi\{B\} = P[W_n \in B]$,  $Q(w;B) = P[W_{n+1}\in B|W_n=w]$, and (the reversibility condition)
\be\label{eq:rvsbl}
	\int_A Q(w;B)\pi\{dw\} = \int_B Q(w;A)\pi\{dw\}
\ee
for Borel sets $A, B \subseteq {\cal W}$, $w \in {\cal W}$, and $n \in \Z$.  Using (and abusing) notation in a standard manner, we write 
$$
	Qf(w) = \int_{\cal W} f(z)Q(w;dz)\ a.e.\ (\pi)
$$ 
for $f \in L^1(\pi)$ and $Q^k = Q\circ \cdots\circ Q$ for the iterates of $Q$.  In addition, let $L_0^{p}(\pi) = \{f \in L^p(\pi): \int_{\cal W} fd\pi = 0\}$,  
$$
	V_n = I+Q+\cdots+Q^{n-1},\quad \bar{V}_n = (V_1+\cdots+V_n)/n,
$$ 
and let $\Vert\cdot\Vert$ denote the norm in an $L^2$ space, either $L^2(\pi)$ or $L^2(P)$.  Finally $\Rightarrow$ denotes convergence in distribution and $\Rightarrow^{\rm p}$ convergence in probability of conditional distributions; that is, if $Z_n:\Omega \to \R$ are random variables and $G$ is a distribution function, then $Z_n|W_0 \Rightarrow^{\rm p} G$, means that the conditional distribution of $Z_n$ given $W_0$ converges in probability to $G$.  

	The reversibility condition (\ref{eq:rvsbl}) is equivalent to requiring $(W_0,W_1)$ and $(W_1,W_0)$ to have the same distribution, since the left side of (\ref{eq:rvsbl}) is $P[W_0 \in A,\ W_1 \in B]$ and the right-hand side is $P[W_0 \in B,\ W_1 \in A]$.  An important consequence (also equivalent) is that the restriction of $Q$ to $L^2(\pi)$ is a self-adjoint operator. For, letting $\langle\cdot,\cdot\rangle$ denote the inner product in $L^2(\pi),\ \langle f,g\rangle = \int_{\cal W} fgd\pi$, $\langle f,Qg\rangle = E[f(W_0)g(W_1)] = E[f(W_1)g(W_0)] = \langle Qf,g\rangle$ for all $f,g \in L^2(\pi)$.
	
	Given $g \in L_0^{2}(\pi)$, let $X_k = g(W_k),\ S_n = X_1+\cdots+X_n$, and $\sigma_n^{2} = E(S_n^{2})$.  Kipnis and Varadhan \cite{KV86} showed that if 
\be\label{eq:var1}
	\lim_{n\to\infty} {\sigma_n^{2}\over n} = \kappa \in [0,\infty), 
\ee
then the conditional distribution of $S_n/\sqrt{n}$ given $W_0$ converges in probability to the normal distribution with mean $0$ and variance $\kappa$.  It is shown in Proposition \ref{prop:var} that $\kappa > 0$ except for trivial special cases;  then $\sigma_n^{-1}S_n|W_0 \Rightarrow^{\rm p} {\rm Normal}[0,1]$. In the proof, Kipnis and Varadhan showed that $S_n$ could be written in the form $S_n = M_n + R_n$, where $0 = M_0,M_1,M_2,\ldots$ is a square integrable martingale with (strictly) stationary increments $D_k = M_k-M_{k-1}$ and $\Vert R_n\Vert = o(\sqrt{n})$.  The result has applications to Markov Chain Monte Carlo, for instance, \cite{T94}, since many algorithms lead to reversible chains; and, to interacting particle systems, \cite{KL99} and \cite{KV86}.  
	
	Here we consider the case in which (\ref{eq:var1}) is weakened to 
\be\label{eq:var2}
	\sigma_n^{2} = n\ell(n),
\ee
where $\ell$ is a slowly varying function, as defined in Chapter 1 of \cite{BGT87}. An example will show that the main result from \cite{KV86} does not extend completely.   Some features do extend, however.  For the remainder of the paper reversibility is assumed along with $g \in L_0^{2}(\pi)$, and $\ell$ is defined by (\ref{eq:var2}).  

	Further developments under the condition (\ref{eq:var1}) may be found in \cite{CP10}; and \cite{MP06} is a recent article on asymptotic normality of sums of stationary processes with non-linear growth of variance.

\section{Generalities}\label{sect:gen}

	In the first proposition, it is shown that only the case $\lim_{n\to\infty} \ell(n) = \infty$ needs to be considered.   The relation
\be\label{eq:var0}
	\sigma_n^{2} = \left[2\langle g,\bar{V}_ng\rangle - \Vert g\Vert^2\right]n 
\ee
is used in its proof.

\begin{proposition}\label{prop:var}
	If $~\liminf_{n\to\infty} \ell(n) < \infty$, then (\ref{eq:var1}) holds; and if $~\liminf_{n\to\infty} \ell(n) = 0$, then $S_n = {1\over 2}[1+(-1)^{n-1}]X_1\ $ with probability one.
\end{proposition}

	{\it Proof}.  Since $Q$ is self-adjoint, we may write $Q = \int_{\Lambda} \lambda dM(\lambda)$, where $\Lambda \subseteq [-1,1]$ is the spectrum of $Q$ and $M$ is a countably additive, projection-valued set function defined on the Borel sets of $\Lambda$.  Then $Q^k = \int_{\Lambda} \lambda^k dM(\lambda)$ for all $k \ge 1$.  See \cite{H54}, Chapter 2.   Let $\mu_g(B) = \langle g,M(B)g\rangle$.  Then $\mu_g$ is a measure, and
\be\label{eq:spct1}
	\langle g,\bar{V}_ng\rangle = \int_{\Lambda} \left[1 - {\lambda\over n}\left({1-\lambda^n\over 1-\lambda}\right)\right] {d\mu_g(\lambda)\over 1-\lambda}.
\ee
Observe that the integrand on the right side of (\ref{eq:spct1}) is non-negative.  So, if $\liminf_{n\to\infty} \ell(n) < \infty$, then the limit inferior of the left side of (\ref{eq:spct1}) is finite and, therefore, 
\be\label{eq:spct2}
	\int_{\Lambda} {\mu_g(d\lambda)\over 1-\lambda} < \infty
\ee
by Fatou's Lemma.  It is clear the integrands on the right side of (\ref{eq:spct1}) are  dominated by an integrable function, hence the integral converges to that on the left side of (\ref{eq:spct2}), and (\ref{eq:var1}) holds with
$$
	\kappa = 2\int_{\Lambda} {d\mu_g(\lambda)\over 1-\lambda} - \Vert g\Vert^2 = \int_{\Lambda} {1+\lambda\over 1-\lambda}\mu_g(d\lambda).
$$
If $\liminf_{n\to\infty} \ell(n) = 0$, then the last integral is $0$ and, therefore, $\mu_g$ is a point mass at $\{-1\}$.  It follows that $Qg = -g$, $E[(X_0+X_1)^2] = 0,\ X_n = (-1)X_{n-1}\ w.p.1$, and $S_n = {1\over 2}[1+(-1)^{n-1}]X_1\ w.p.1$. 

	As a consequence there is no loss of generality in supposing that $\ell(n) \to \infty$, and we shall do so where convenient.  For if $\liminf_{n\to\infty} \ell(n) < \infty$, then the Kipnis-Varadhan result is applicable.

	The proof of the next proposition uses (\ref{eq:var0}) and
\be\label{eq:vbar}
	\bar{V}_n = \sum_{k=0}^{n-1} \left(1-{k\over n}\right)Q^k.
\ee

\begin{proposition}\label{prop:vee}
	If $\ell$ varies slowly  in (\ref{eq:var2}), then $\Vert V_ng\Vert = o(\sigma_n)$.
\end{proposition}

	{\it Proof}. Using the reversibility and (\ref{eq:vbar}),
\begin{eqnarray*}
	 \Vert V_n g\Vert_2^{2} &=& \sum_{j=0}^{n-1}\sum_{k=0}^{n-1} \< g, Q^{k+j} g \> \\
      		&=& \sum_{i=0}^{n-1}(i+1)\<g,Q^i g \>+\sum_{i=n}^{2n-2}(2n-1-i)\< g,Q^i g\>  \\
   		&=& \sum_{i=0}^{2n-2} (2n-1-i) \< g,Q^i g \> - 2 \sum_{i=0}^{n-1}(n-1-i) \< g, Q^i g \>\\
    		&=& {1\over 2} \left[ \sigma_{2n-1}^2+(2n-1) \Vert g\Vert^2 \right]-\left[\sigma_{n-1}^2+(n-1)\Vert g\Vert^2\right]\\
    	&=& {1\over 2}\sigma_{2n-1}^2-\sigma_{n-1}^2 +{1\over 2} \Vert g \Vert^2.
\end{eqnarray*}
The proposition then follows directly from (\ref{eq:var2}) and the slow variation of $\ell$.

\begin{corollary}
	If $\ell$ varies slowly, then there is a sequence of square integrable martingales $0 = M_{n,1},M_{n,2},\ldots$ with stationary increments $D_{n,k} = M_{n,k}-M_{n,k-1}, k\ge 1,$ for which $\max_{k\le n} \Vert S_k-M_{n,k}\Vert = o(\sigma_n)$.
\end{corollary}

	{\it Proof}.  This follows from Proposition \ref{prop:vee} and Theorem 1 of \cite{WW04}. It is relevant that 
$$
	D_{n,k} = \bar{V}_ng(W_k)-Q\bar{V}_ng(W_{k-1})
$$ 
and $M_{n,k} = D_{n,1}+\cdots+D_{n,k}$ in the proof of the latter result.

\begin{corollary}\label{cor:mrtclt}
	If $\ell$ varies slowly and there is a $\lambda \ge 0$ for which
\be\label{eq:stbl}
	{1\over\sigma_n^{2}}\sum_{k=1}^n E(D_{n,k}^2|W_{k-1}) \to^{\rm p}\lambda
\ee
and
\be\label{eq:lfc}
	{1\over\sigma_n^{2}}\sum_{k=1}^n E(D_{n,k}^2{\bf 1}_{\{|D_{n,k}|>\epsilon\sigma_n\}}|W_{k-1}) \to^{\rm p} 0
\ee
for every $\epsilon > 0$, then 
\be\label{eq:an}
	{S_n\over\sigma_n}|W_0 \Rightarrow^{\rm p} {\rm Normal}[0,\lambda].
\ee
\end{corollary}

	{\it Proof}.  This follows from the Martingale Central Limit Theorem, e.g. \cite{B95}, pp. 475-478, applied conditionally given ${\cal F}_0:=\sigma(\ldots W_{-1}, W_0)$.  For $\lambda = 1$ the proof is detailed in \cite{WW04}, and the extension to $\lambda \ne 1$ presents no difficulty. 

	In the next proposition we write $S_n = S_n(g)$ and $\sigma_n = \sigma_n(g)$ to emphasize the dependence on $g$.	  We also use the following:

\begin{lemma}\label{lem:sltzky}
	If $Z_n|W_0 \Rightarrow^{\rm p} G$ and $Z_n'-Z_n \to^{\rm p} 0$, then $Z_n'|W_0 \Rightarrow^{\rm p} G$.
\end{lemma}

	{\it Proof}.  Lemma \ref{lem:sltzky} follows from the unconditional version of Slutzky's Theorem, by considering subsequences along which convergence in probability can be replaced by almost sure convergence.

\begin{proposition}\label{prop:rigid}
	If $\ell(n) \to \infty$, and (\ref{eq:an}) holds for a given $g$, then for any $j \ge 1$, $\sigma_n(Q^jg) \sim \sigma_n(g)$ and (\ref{eq:an}) holds with the same $\lambda$ when $g$ is replaced by $Q^jg$.
\end{proposition}

	{\it Proof}.  It suffices to prove the corollary for $j = 1$; and in this case it follows from $S_n(g)-S_n(Qg) = \sum_{k=1}^n [g(W_k)-Qg(W_{k-1})] + Qg(W_0)-Qg(W_n)$, which implies
$$
	\eqalign{|\sigma_n(g)-\sigma_n(Qg)| &\le \Vert S_n(g)-S_n(Qg)\Vert\cr 
			&\le \Vert g(W_1)-Qg(W_0)\Vert\sqrt{n} + 2\Vert Qg(W_0)\Vert = o[\sigma_n(g)],\cr}
$$
and Lemma \ref{lem:sltzky} above. 

\medskip\noindent
{\bf Remark 1}.  The proof of the Proposition \ref{prop:rigid} did not use the reversibility and, therefore, is valid for any stationary process.
	
\medskip \noindent
	{\bf Remark 2}.  Proposition \ref{prop:rigid} illustrates an important difference between the case $\ell(n) \to \infty$ and $\ell(n) \to \kappa$, considered in \cite{KV86}.
For if (\ref{eq:var1}) holds, then
\begin{equation} \label{eq:kappa}
	\kappa = \kappa(g) = 2 \lim_{n\to\infty} {\sum_{k=0}^{n}} \left(1-{k\over n}\right)\langle g,Q^kg\rangle - \Vert g\Vert^2,
\end{equation}
It is then not difficult to see that (\ref{eq:kappa}) holds  when $g$ is replaced by  $Q^jg$; and $[\kappa(g)+\cdots+\kappa(Q^ng)]/n$ approaches zero as $n \to \infty$, by Theorem {2} of \cite{ZW08}.
	
	\medskip\noindent
	{\bf Remark 3}.   Kipnis and Varadhan showed that if (\ref{eq:var1}) holds then $D_{n,k}$ converges in $L^2(P)$ for every $k$.  Clearly, this is impossible if $\ell(n) \to \infty$.  If it were the case, however, that $D_{n,1}/\sqrt{\ell(n)}$ converged in $L^2(P)$, then (\ref{eq:stbl}) and (\ref{eq:lfc}) would follow easily with $\lambda = 1$, and the conditional distributions of $S_n/\sigma_n$ would converge to the standard normal distribution, as noted in \cite{WW04}.  This hope cannot be realized either, however, if $\lim_{n\to\infty} \ell(n) = \infty$.  For, $D_{n,1}/\sqrt{\ell(n)}$ cannot be a Cauchy sequence, in this case.  To see this first observe that 
$$
		\left\Vert \frac{D_{n,1}}{\sqrt{\ell(n)}}-\frac{D_{m,1}}{\sqrt{\ell(m)}} \right\Vert^{2}=\frac{1}{\ell(n)} \Vert D_{n,1} \Vert^2+\frac{1}{\ell(m)}\Vert D_{m,1}\Vert^2-\frac{2}{\sqrt{\ell(m)\ell(n)}} \<D_{m,1},D_{n,1} \>
$$
and
\begin{eqnarray*}
   \<D_{m,1}, D_{n,1} \> &=& \<\bar{V}_n g(w_1)-Q\bar{V}_n g(w_0), \bar{V}_m g(w_1)-Q \bar{V}_m g(w_0) \> \\
   &=& \<\bar{V}_n g, \bar{V}_m g \>-\< Q\bar{V}_n g, Q \bar{V}_m g \>  \\
   &=& \<\bar{V}_n g, \bar{V}_m g \>-\<Q^2 \bar{V}_n g, \bar{V}_m g \> \\
   &=& \<(I-Q^2) \bar{V}_n g, \bar{V}_m g \>\\
   &=& \<(V_2-\frac{1}{n} Q V_n V_2)g, \bar{V}_m g \>.
\end{eqnarray*}
So, for any fixed $m$,
$$
  \lim_{n\to \infty}  \left\Vert \frac{D_{n,1}}{\sqrt{\ell(n)}}-\frac{D_{m,1}}{\sqrt{\ell(m)}} \right\Vert^{2} =1+ \frac{1}{\ell(m)}\Vert D_{m,1} \Vert^{2},
$$
and, therefore,
$$
   \lim_{m\to \infty}\lim_{n\to \infty}  \left\Vert \frac{D_{n,1}}{\sqrt{\ell(n)}}-\frac{D_{m,1}}{\sqrt{\ell(m)}} \right\Vert^{2}  = 2.
$$

\section{Examples}\label{sect:ex}

		For a simple reversible chain, let $\nu$ be a probability measure on the Borel sets of $\R$ and $p:\R \to (0,1)$ a measurable function for which
$$
	\theta = \int_{\R} {d\nu\over 1-p} < \infty;
$$
and let
\be\label{eq:clsscl}
	Q(w;B) = p(w){\bf 1}_B(w) + [1-p(w)]\nu\{B\}
\ee
for Borel set $B \subseteq \R$ and $w \in \R$.  Then $Q$ is a stochastic transition function with stationary distribution 
$$
	d\pi = {d\nu\over\theta(1-p)},
$$
and (\ref{eq:rvsbl}) is satisfied.  Thus there is a reversible Markov Chain $\ldots W_{-1},W_0,W_1,\ldots$ with transition function $Q$ and marginal distribution $\pi$.  This construction is classical and is described in \cite{R00}, pp. 134-135.

	Now let $\tau_0,\tau_1,\tau_2,\ldots$ be the times before the process jumps, $\tau_0 = \max\{n \ge 0: W_n = W_0\}$ and 
$$
	\tau_{k} = \max\{n > \tau_{k-1}: W_n = W_{\tau_{k-1}+1}\}.
$$
Then $W_{\tau_k} = W_{\tau_{k-1}+1}$, and
$$
	S_{\tau_m} = \tau_0X_0 + (\tau_1-\tau_0)X_{\tau_1}+\cdots+(\tau_m-\tau_{m-1})X_{\tau_m}.
$$
By the Markov property, $(\tau_0,W_0)$ and $[(\tau_j-\tau_{j-1}),W_{\tau_j}],\ j \ge 1$ are independent random vectors for which $W_{\tau_j} \sim \nu $ and
$$
	P[\tau_j-\tau_{j-1} \ge k|W_{\tau_j}=w] = p(w)^{k-1}
$$
for all $w \in {\cal W}$, $k \ge 1$ and $j \ge 1$.  It follows that $E(\tau_j-\tau_{j-1}|W_{\tau_j}=w) = 1/[1-p(w)]$ and 
$$
	E(\tau_j-\tau_{j-1}) = \int_{\cal W} {d \nu \over1-p} = \theta.
$$
By way of contrast, $W_{\tau_0} = W_0 \sim \pi$, and $E(\tau_0) = \int pd\pi/(1-p)$, possibly infinite.  Let $Y_j = (\tau_j-\tau_{j-1})X_{\tau_j}$ and $T_m = Y_1+\cdots+Y_m$, so that $S_{\tau_m} = \tau_0W_0 + T_m$.   Then $Y_1,Y_2,\ldots$ are independent and identically distributed; moreover, $E(Y_j) = 0$, since 
$$
	E\left[(\tau_j-\tau_{j-1})X_{\tau_j}\right] = E\left[{g(W_{\tau_j})\over 1-p(W_{\tau_j})}\right] = \int_{\cal W} {g\over 1-p}d \nu  = \theta\int_{\cal W} gd\pi,
$$
  and $g \in L_0^{2}(\pi)$.  Let 
$$
	 H(y) = \int_{|Y_j|\le y} Y_j^{2}dP ,
$$
and recall the following version of the Central Limit Theorem for i.i.d. variables (with possibly infinite variances), for example, \cite[pp. 576-578]{F71}: If $Y_1,Y_2,\ldots$ are (any) i.i.d. random variables for which $E(Y_j) = 0$ and $ H (y)$ varies slowly at $\infty$, then there are $\gamma_m$ for which 
$$
	\gamma_m^{2} \sim m H(\gamma_m)\quad {\rm and}\quad {T_m\over \gamma_m} \Rightarrow {\rm Normal}[0,1].
$$  

	The following lemma is intuitive.  The proof is presented after Proposition \ref{prop:clt} is established. To state it,  define integer-valued random variables $m_n$ such that  $\tau_{m_n} \le n < \tau_{m_{n}+1}$ for $n=1,2\ldots$ .

\begin{lemma}\label{lem:clt}
	As $n \to \infty,\ S_n - T_{m_n} = O_p(1)$; and if $ H $ varies slowly at $\infty$, then $T_{m_n}-T_{\lfloor n/\theta\rfloor} = o_p(\gamma_n)$.
\end{lemma}
	
\begin{proposition}\label{prop:clt}
	If $ H (y)$ varies slowly and $\gamma_m^{2} \sim m H (\gamma_m)$, then $$
		{S_n\over\gamma_n}|W_0 \Rightarrow^{\rm p} {\rm Normal}\left[0,{1\over\theta}\right].
$$
\end{proposition}

	{\it Proof of Proposition \ref{prop:clt}}. \quad That $T_m/\gamma_m \Rightarrow {\rm Normal}[0,1]$ was noted above.  So, since $\gamma_{\lfloor n/\theta\rfloor} \sim \gamma_n/\sqrt{\theta}$, $T_{\lfloor n/\theta\rfloor}/\gamma_n \Rightarrow {\rm Normal}[0,1/\theta]$; and since $W_0$ and $T_m$ are independent for all $m$, the conditional distributions have the same limit.  The proposition now follows directly from Lemmas \ref{lem:sltzky} and \ref{lem:clt}. 
	
		{\it Proof of Lemma \ref{lem:clt}}.\quad First observe that $S_n -  T_{m_n} = \tau_0W_0 + (n-\tau_{m_n})W_{\tau_m+1}$.  It is clear that $\tau_0W_0$ is stochastically bounded and that $|(n-\tau_{m_n})W_{\tau_{m_n+1}}| \le (\tau_{m_n+1}-\tau_{m_n})|W_{\tau_{m_n+1}}|$. To see that the latter term is stochastically bounded, let $f$ denote the marginal mass function of $\tau_j-\tau_{j-1},\ j \ge 1$.  Then the asymptotic distribution of $\tau_{m_n+1}-\tau_{m_n}$ has probability mass function $\tilde{f}(k) =kf(k)/\theta$, by the Renewal Theorem, \cite[p.271]{F71}, and the conditional distribution of $W_{\tau_{m_n+1}}$ given $\tau_{m_n+1}-\tau_{m_n}$ does not depend on $n$.  That $(\tau_{m_n+1}-\tau_{m_n})|W_{\tau_{m_n+1}}| = O_p(1)$ follows easily.

	The proof of the second assertion uses the following version of L\' evy's Inequality \cite[p.259]{L63}: If $H$ varies slowly at $\infty$, then $K^{-1} := \inf\{\min[P(T_k<0),P(T_k>0)]: k \ge 1\} > 0$, and
\be\label{eq:levy}
	P\left[\max_{k\le n} |T_k| > t \right] \le KP[|T_n| > t]
\ee
for all $t > 0$.  Observe that
\be\label{eq:anscmb}
	\eqalign{P\left[|T_{m_n}-T_{\lfloor {n\over\theta}\rfloor} | \ge \epsilon\gamma_n \right] &\le P\left[|m_n-\lfloor {n\over\theta}\rfloor| \ge \delta n\right]\cr
		&+ P\left[ \max_{|k\theta -n| \le \theta\delta n+\theta} |T_k - T_{\lfloor {n\over\theta}\rfloor}| \ge \epsilon\gamma_n\right] .\cr}
\ee
The first term on the right approaches $0$ for any $\delta > 0$ by the Law of Large Numbers.  Letting $N_n = \lfloor n\delta/\theta\rfloor + 4$ and using (\ref{eq:levy}), the second is at most $2KP[|T_{N_n}| \ge \epsilon\gamma_n]$.  So, by the Central Limit Theorem, the limit superior of the right side of (\ref{eq:anscmb}) is at most $4K[1 - \Phi({\epsilon/\sqrt{\delta}})]$, which approaches $0$ as $\delta \to 0$. 

	For the example below, observe that if $f \in L^1(\pi)$, then $Qf(w) = p(w)f(w) + [1-p(w)]\int_{\cal W} fd \nu$.  So, if ${\cal W} = \R$, $ \nu $ is a symmetric measure, $p$ is a symmetric function, and $f$ is an odd function, then $Q^nf = p^n\times f$.

\begin{ex}\label{ex:nrml}{\rm
	Consider (\ref{eq:clsscl}) with, $p(w) = e^{-{1/|w|}}$, 
\be\label{eq:psi}
	\nu \{dz\} = e{[1-p(z)]dz\over 2z^2}\ {\rm for}\ |z| \ge 1,
\ee
in which case $\theta = e$ and $\pi\{dw\} = {dw/2w^2}\ {\rm for}\ |w| \ge1$.  Let $g(w) = {\rm sign}(w)$.  Then $g \in L_0^{\infty}(\pi)$, $Q^ng = p^n\times g$, and 
$$
	\langle g,Q^ng\rangle = \int_{\R} p^nd\pi = \int_{1}^{\infty} e^{-{n\over w}}{dw\over w^2} = {1\over n}\int_{0}^n e^{-x}dx \sim {1\over n}.
$$
It follows that $\langle g,\bar{V}_ng\rangle \sim \langle g,V_ng\rangle \sim \log(n)$ and 
$\sigma_n^{2} = [2\langle g,\bar{V}_ng\rangle -\Vert g\Vert^2]n \sim 2 n\log(n)$.  So, (\ref{eq:var2}) is satisfied with $\ell(n) \sim 2\log(n)$.

	Recall the definition of the $\tau_j$ and the distribution of $[\tau_j-\tau_{j-1},W_{\tau_j}]$.  Then 
$$
	P[|(\tau_j-\tau_{j-1})X_{\tau_j}| > k] = P[(\tau_j-\tau_{j-1}) > k] = \int_{\R} p^{k}d \nu  = e\int_{\R} (1-p)p^kd\pi.
$$
The last integral in the previous display is just
$$
	\int_{1}^{\infty} \left(1-e^{-{1\over z}}\right) e^{-{k\over z}}{dz\over z^2} = {1\over k}\int_{0}^k \left(1-e^{-{y\over k}}\right) e^{-y}dy \sim {1\over k^2}\int_{0}^{\infty} ye^{-y}dy = {1\over k^2};
$$
thus,
\be\label{eq:tail}
	P[|(\tau_j-\tau_{j-1})X_{\tau_j}| \ge k] \sim {e\over k^2}.
\ee
It follows easily that $ H(y) \sim 2e\log(y) = e\ell(y)$, $\gamma_n^{2} = 2en\log(\gamma_n) \sim en\log(n) = {1\over 2}e\sigma_n^{2}$, and, therefore,
$$
	{S_n\over\sigma_n} \Rightarrow {\rm Normal}\left[0,{1\over 2} \right]
$$
(a non-standard normal distribution). 

	Since $E(\sigma_n^{-2}S_n^2)$ is bounded, it follows that $E|S_n| \sim \pi^{-{1\over 2}}\sigma_n$ and, therefore, that $S_n/E|S_n| \Rightarrow {\rm Normal}[0,{1\over 2}\pi]$.  The latter convergence can also be deduced from Theorem 4 of \cite{MP06}.  To do so, it suffices to verify Equation (3.2) of that paper: Since $|g| \le 1$, it is not difficult to see that the term whose limit is taken in (3.2) is at most $\sigma_n^{-2}\sum_{k=1}^n k\beta_k$, where $\beta_k$ is the coefficient of absolute regularity.  So, it suffices to show that $\beta_n$ is of order $1/n$, and this may be deduced from the equation at the top of page 136 of \cite {R00} together with the relation $P[\tau_0 > n] = \int_{\R} p^nd\pi \sim 1/n$.  (The $\tau$ in \cite{R00} is our $\tau_0+1$.)  Conditional convergence is not asserted in Theorem 4 of \cite{MP06} but is implicit in the proof; that $E|S_n| \sim \pi^{-{1\over 2}}\sigma_n$ is not deducible from that theorem, however, because $S_n$ is not normalized by $\sigma_n$ there.

}

\end{ex} 

\begin{ex} {\rm A slight modification of Example 1 produces a very simple bounded stationary sequence whose normalized partial sums converge in distribution to a stable distribution. Other examples may be found in \cite{KS01}.  If (\ref{eq:psi}) is changed to
$$
	\nu \{dz\} = {[1-p(z)]dz\over 2\gamma_{\alpha} |z|^{\alpha}}
$$
for $|z| \ge 1$, where $1 < \alpha < 2$ and $\gamma_{\alpha} = \int_{0}^{1} y^{\alpha-2}(1-e^{-y})dy$, then $\pi\{dz\} = (\alpha-1)/(2 |z|^{\alpha})dz$ for $|z| \ge 1$, and 
$$
	P\left[Y > y\right] \sim {\Gamma(\alpha)\over\gamma_{\alpha}y^{\alpha}}
$$
as $y \to \infty$.  It then follows that $n^{-{1\over\alpha}}S_n \Rightarrow Z$, where $Z$ has a symmetric stable distribution with characteristic function $e^{-c_{\alpha}|t|^{\alpha}}$ and $c_{\alpha} = (\alpha-1)\Gamma(\alpha)\int_{0}^{\infty} x^{-\alpha}\sin(x)dx$.  
}
\end{ex}

\end{document}